\documentclass{amsart}
\usepackage{amsthm}
\theoremstyle{plain}
\newtheorem{theorem}{Theorem}[section]

\newtheorem{fact}[theorem]{Fact}

\newtheorem{corollary}[theorem]{Corollary}
\theoremstyle{definition}
\newtheorem{definition}[theorem]{Definition}

\usepackage[T1]{fontenc}
\usepackage{graphicx}
\usepackage{hyperref}
\usepackage{color}

\usepackage{amsmath}
\usepackage{amssymb}
\usepackage{mathtools}
\usepackage{extpfeil}
\usepackage{cite} 
\usepackage{tikz}
\usetikzlibrary{automata, positioning, cd, shapes.geometric, backgrounds}

\DeclareMathOperator{\image}{Im}
\newcommand{\tinysupset}{\mathrel{\scalebox{.4}[.4]{$\supset$}}}
\newcommand{\cp}{\wr_{\tinysupset}}

\newcommand\dhxrightarrow[2][]{%
  \mathrel{\ooalign{$\xrightarrow[#1\mkern4mu]{#2\mkern4mu}$\cr%
  \hidewidth$\rightarrow\mkern4mu$}}
}

\begin{document}


\title[From Relation to Emulation and Interpretation]{From Relation to
Emulation and Interpretation: Computer Algebra Implementation of the Covering Lemma for Finite Transformation Semigroups}
\author{Attila Egri-Nagy$^1$}
\address{$^1$Akita International University, Japan}
\email{egri-nagy@aiu.ac.jp}
\author{Chrystopher L. Nehaniv$^2$}
\address{$^2$University of Waterloo, Canada}
\email{cnehaniv@uwaterloo.ca}

\maketitle              

\begin{abstract}

We give a practical computer algebra implementation of the Covering Lemma for finite
transformation semigroups.
The lemma states that given a surjective relational morphism
$(X,S)\twoheadrightarrow(Y,T)$,
we can establish emulation by  a cascade
product (subsemigroup of the wreath product): $(X,S)\hookrightarrow (Y,T)\cp (Z,U)$.
The dependent component $(Z,U)$ contains the kernel of the morphism, the information lost in the map.

The implementation complements the existing tools for the holonomy decomposition
algorithm.
It gives an incremental method to get a coarser decomposition when
computing the complete skeleton for holonomy is not feasible.
Here, we describe a simplified and generalized algorithm for the lemma and compare it to the
holonomy method.
Incidentally, the kernel-based method could be the easiest way of
understanding the hierarchical decompositions of transformation semigroups and
thus the celebrated Krohn-Rhodes theory.
\keywords{transformation semigroup \and relational morphism \and cascade product.}
\end{abstract}

\section{Introduction}

The Krohn-Rhodes theory \cite{primedecomp65} is about how to decompose automata and build them from elementary components.
As a fundamental theory, its applicability goes beyond algebra and automata theory \cite{wildbook}.
For applications, we need computational tools.
The \texttt{SgpDec} \cite{SgpDec2014} package for the \texttt{GAP} system \cite{GAP4} provides an implementation for the holonomy decomposition method \cite{holonomy_algorithm}.
This computer algebra tool enabled the investigation of several applications of
the decompositions (see, e.g., \cite{PhilTransA2015}).
However, the limitations of the holonomy algorithm also became clear. It has exponential resource usage: the number of images of the state set is bounded by $2^n$ for $n$ states.
At this price, it gives a `highest resolution' decomposition, which is desired in some situations, but not in others.
Hence, there is a need for a more flexible and incremental decomposition algorithm.

Here, we use an alternative proof of the Krohn-Rhodes Theorem, the Covering Lemma from \cite{FromRelToEmulationNehaniv1996} and simplify it for a computer algebra implementation.
Starting with a surjective homomorphism, we build a two-level decomposition.
The top level component contains the information preserved by the morphism, and the bottom component the information left out.
The latter is the \emph{kernel} of a homomorphism, which is, in general, the inverse image of the identity.
For groups, it is a subgroup.
For semigroups, we have a more general, possibly partial algebraic structure, a category or a semigroupoid.
This kernel, called the `derived category', is an important tool in advanced semigroup theory \cite{Tilson1987}.
In \cite{FromRelToEmulationNehaniv1996}, the kernel maps functorially to sets
and characterizes the minimal computation needed to ``undo'' a relational morphism.
Here, we take a more elementary approach, omitting the category theoretical concepts.

\section{An algebraic view of computation -- Definitions}
We model finite state transition structures (computers) algebraically, as
transformation semigroups \cite{compstructs2017}.
We focus on \emph{how} the computation is done, the dynamics, instead of only its result.
We do not specify initial and accepting states, thus to avoid confusion, we do not call transformation semigroups automata.
The following definitions will describe what is a computational structure, how to build more complex ones out of simpler components, and how to compare their computational capabilities.

\subsection{Transformation and Cascade Semigroups}

A \emph{transformation semigroup} $(X,S)$ is a finite nonempty set of
\emph{states} $X$  and a set $S$ of total transformations of $X$, i.e.,
functions of type $X\rightarrow X$, closed under composition.
We use the symbol $\cdot$ both for the semigroup action $X\times S\to X$, e.g.,
$x\cdot s$, and for the semigroup multiplication $S\times S\to S$, e.g., $s\cdot
t$.
We also extended the operations for sets $A\cdot B=\{ a\cdot b \mid a\in A, b \in B\}$.
Alternatively, we omit the symbol, so juxtaposition is acting on the right when
computing with transformations.
For representing states we use positive integers.
We denote constant maps as $c_j$, where $j$ is the unique point in its image set.
The action is \emph{faithful} if $(\forall x\in X, xs=xs')\implies
s=s'$. 

We say that $A\subseteq S$ generates $S$, denoted by $\langle A \rangle=S$, if
$S$ is the smallest semigroup containing $A$.
In other words, $S$ consists of the products of elements of $A$.
In automata terminology, the generators of the semigroup correspond to the input
symbols.
For efficient computations, we aim to work with a small set of generators
only, avoiding the complete enumeration of semigroup elements.

%

\subsection{Cascade Products}

The {\em wreath product} $(X,S)\wr (Y,T)$  of transformation semigroups
is the \emph{cascade transformation semigroup} $(X\times Y,W)$ where
$$W=\{(s,d) \mid  s\in S, d\in T^X\},$$
whose elements map  $X\times Y$ to itself as follows
$$(x,y)\cdot (s,d)=(x\cdot s, y\cdot d(x))$$ for $x\in X, y\in Y$.
We call $d: X\to T$ a dependency function.
Note that it maps states to transformations.
It shows how the action in the bottom component depends on the top level state.
The term `cascade' refers to the flow of control information from the top level component to the bottom level component: the state $x$ determines the action in $T$.
Here, $T^X$ is the semigroup of  all functions  $d$ from $X$ to $T$, thus the wreath product has exponential size, and it is unsuitable for computer implementations.

We call a (proper) subsemigroup of the wreath product a \emph{cascade product}.
We denote such a product operation by $\cp$, indicating that it is a substructure and also hinting the direction of the control flow:
$$ (X,S) \wr (Y,T) \supset (X,S)\cp (Y,T).$$
The wreath product is the full cascade product, and it is uniquely determined by its components.
On the other hand, there are cascade products as many as subsemigroups of the wreath product.
Thus, the cascade product notation could be ambiguous, and the context should make it clear which one are we talking about.
In practice, we specify a cascade product by the set of transformations, or by a
generator set, or as an image of a relation.

\subsection{Relational Morphisms}
Here we define algebraic structure-preserving maps, but in terms of relations, not functions.
In semigroup theory there are several reasons to do this \cite{QBook,MFAT2020}.
Since we work with transformation semigroups, we need to define two relations:
one for states and the other for transformations.
\begin{definition}[Relational Morphism]
  A \emph{relational morphism} of transformation semigroups $(X,S)\xrightarrow{\theta,\varphi} (Y,T)$ is a pair of relations $(\theta: X\rightarrow Y, \varphi: S \rightarrow T)$ that are fully defined, i.e., $\theta(x)\neq\varnothing$ and $\varphi(s)\neq\varnothing$, and satisfy the condition of compatible actions for all $x\in X$ and $s\in S$:
$$y\in\theta(x), t\in \varphi(s) \implies y\cdot t \in \theta(x\cdot s),$$
or more succinctly: $\theta(x)\cdot\varphi(s)\subseteq \theta(x\cdot s).$
\end{definition}
We will generally be assuming that $\varphi$ is a relational morphism of
semigroups, i.e., $\varphi(s)\cdot\varphi(s')\subseteq \varphi(s s')$. 
The relations can also be expressed as set-valued functions: 
$\theta : X\rightarrow
2^Y$, $\varphi : S\rightarrow 2^T$.
In this paper we refer to relational morphisms simply as morphisms, and to avoid misunderstanding we never abbreviate the word `homomorphism'.

A morphism is \emph{surjective} if $\bigcup_{x\in X}\theta(x)=Y$ and
$\bigcup_{s\in S}\varphi(s)=T$.
In computational implementations surjectivity is often given, as the image
transformation semigroup is defined by the morphism.
A surjective morphism can lose information, thus we call the target an \emph{approximation} of the source.
A morphism is \emph{injective} if $\theta(x)\cap\theta(x')\neq \varnothing$ implies $x=x'$, and
analogously $\varphi(s)\cap\varphi(s')\neq \varnothing$ implies $s=s'$.
In other words, the relations map to disjoint sets.
Injective morphism is the same as a \emph{division}.
Injectivity guarantees that no information is lost during the map, therefore the
target of the morphism can compute at least as much as the source.
Thus, we call injective morphisms \emph{emulation} or \emph{covering}.
We use the standard arrow notations, $\hookrightarrow$ for injective, and $\twoheadrightarrow$ for surjective morphisms.

For relations, being surjective and injective at the same time does not imply
isomorphism, but it implies being the inverse of a homomorphism.
On the other hand, we have the ability to turn around surjective morphisms to get an injective ones and vice-versa.
Relational morphisms preserve subsemigroups both directions.

\section{The General Algorithm}
\label{sect:generalalg}
First we describe the general mechanism of how parts of the semigroup action in $(X,S)$ are transferred into the newly constructed $(Z,U)$.

\subsection{Problem Description}
Our goal is to understand the inner workings of a transformation semigroup $(X,S)$.
We would like to build a hierarchical decomposition, a cascade product that emulates $(X,S)$.
We get simpler components by distributing the computation over the levels.
In general, our task is to find the unknowns in
$$(X,S)\hookrightarrow (?,?)\cp(?,?)\cp\ldots \cp (?,?).$$
Requiring only the emulation leaves the solution space huge.
The holonomy decomposition uses constraints from the detailed examination of how $S$ acts on the subsets of $X$ and produces a more limited set of solutions \cite{holonomy_algorithm}.
They are essentially the same up to choosing equivalence class representatives and moving the components around when their dependencies allow some wiggle room.

Here, we take a more flexible approach.
We define the top level component, 
$$(X,S)\hookrightarrow (Y,T)\cp(?,?),$$
and leave it to the algorithm to define the bottom level component, recovering all the information about $(X,S)$ that is not represented in $(Y,T)$.
The decomposition has two levels, but the process can be iterated to have a more fine-grained solution.

\subsection{Constructing the Emulation}

Given a relational morphism as an input we build an emulation as output:
$$\text{from } \ (X,S)\dhxrightarrow{R(\theta,\varphi)} (Y,T) \ \text{ to } \ (X,S)\xhookrightarrow{E(\psi,\mu)}(Y,T)\cp(Z,U).$$
We need to construct the transformation semigroup $(Z,U)$, the dependency functions $Y\to U$ for defining the cascade transformations.
In practice, we will define only the emulation relational morphism $E(\psi, \mu)$ and only map the generators to define the cascade product.

\emph{What are the states in $Z$?}
We need to specify one or more coordinate pairs for $x\in X$.
The top level coordinates are elements of $\theta(x)\subseteq Y$.
Each such top level state $y$ gives a context for
the bottom level states.
In each context, what a particular $z\in Z$ refers to is given by a \emph{labelling}
mechanism, a partial one-to-one mapping $w_y:X\to Z$, that depends on $y$.
More precisely, we only need to map $\theta^{-1}(y)\subseteq X$ to $Z$.
The states in $Z$ can be reused in the different contexts,  we need at least $|Z|=\max_y|\theta^{-1}(y)|$.
The coordinate pair $(y,z)$ gives the context $y$ and the encoded state $z$ as
lift for $x$.
Formally,
$$ \psi(x)=\{ (y,z) \mid y \in \theta(x), z=x\cdot w_y \}.$$

\emph{How to construct transformations in $U$?}
Using the labelling above, we simply {\em transfer} the relevant bits of $s$ acting on
$X$:
Formally, to see where a single state $z$ should go respecting the original
action of $s\in S$, we do
 $z\mapsto z\cdot w_y^{-1}\cdot s \cdot w_{yt}$,
where $w_y^{-1}:Z \rightarrow \theta^{-1}(y)$ is a partial function, $s:X\to X$, and
$w_y:\theta^{-1}(y)\to Z$.
The net result of this composition is a map $Z\to Z$.
If we want to identify the main trick of the method, then this is it: decode the
states according to the current context, perform the original action, then
encode the states for the new context.
We have a `window' through $\theta^{-1}(y)$ into $X$, or using another metaphor,
$\theta^{-1}(y)$ sheds light on some part of $X$, so we can see the action of
$S$ there.
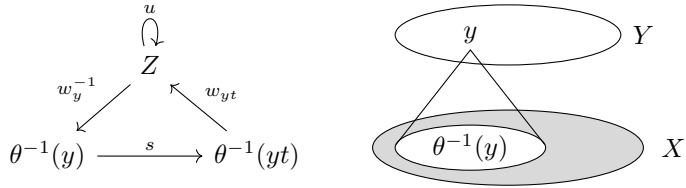
\begin{center}
   \begin{tikzcd}[column sep=small,background grid/.style={thick,draw=gray!20,step=.5cm},
    ]
   & Z \arrow[dl,"w_y^{-1}"'] \arrow[loop above, "u"] & \\
\theta^{-1}(y) \arrow[rr, "s"] & & \theta^{-1}(yt) \arrow[ul, "w_{yt}"']
\end{tikzcd}
\begin{minipage}{.5\textwidth}
\begin{tikzpicture}
[background grid/.style={thick,draw=gray!20,step=.5cm},
]
\draw [use as bounding box,draw=none] (0,0) rectangle (5,2.5);
   \draw [draw](2.5,2.2) ellipse (1.5 and .4);
   \draw [draw, fill=gray!30](2.5,0.7) ellipse (1.8 and .5);
   \draw [draw, fill=white](2,0.7) ellipse (1 and .3);
   \node at (2,2.2) {$y$};
   \node at (2,0.7) {$\theta^{-1}(y)$};
   \node at (4.3,2.2) {$Y$};
   \node at (4.7,0.7) {$X$};
\draw [draw] (1.01,0.74) -- (2,2.0) -- (2.99,0.74);
 \end{tikzpicture}
 \end{minipage}
\end{center}
\emph{How is it possible at all that differently encoded local actions compose nicely when multiplying cascade transformations?}
Firstly, in the same context, for a fixed $y$, they compose naturally. For each
$y$ we have a differently behaving clone of the component (a view advocated in
\cite{bitvectorcascade}), but as long as the top level is fixed, we stay in the
same clone transformation semigroup.
Secondly, the context switching is built-in in $u=w_y^{-1} s w_{yt}$, since the result of $s$ is encoded for the $yt$ context.
In short, when lifting $s$ into the cascade product we have
$$s\mapsto\{ (t, w_y^{-1} s w_{yt} ) \mid t\in \varphi(s), y\in \image\theta\}.$$
Now, if we consider a sequence of actions, it is easy to see how the (partial) labelling function in the middle matches its inverse, thus cancelling out.
$$ss'\mapsto\{ (tt', w_y^{-1} s w_{yt} w_{yt}^{-1}s'w_{ytt'}=w_y^{-1} ss'w_{ytt'} ) \mid t\in \varphi(s), t'\in\varphi(s'),  y\in \image\theta\}.$$

To see how this forms a valid cascade transformation, we can see that $t\in T$ by the definition of $\varphi$. 
The bottom level action  $w_y^{-1} s w_{yt}$ has the form of a dependency function $Y\to U$, since for an $(s,t)$ pair the expression's value depends on $y$, and $w_y^{-1} s w_{yt}$ as a function has the type $Z\to Z$ by the definition of the labelling.
Then, we define lifting to a cascade by
$$\mu(s)=\{(t, w_y^{-1} s w_{yt}) \mid t\in \varphi(s), y\in \image\theta\}.$$
Therefore, we will have a cascade transformation for each image $t$ of $s$ in $\varphi$, and each cascade will be defined for every possible state image $y$ by $\theta$. In practice, we simply have $Y=\image\theta$.

To summarize, \emph{how does this method recover the information lost through $R$?} For collapsed states, we go back to
the preimages of $y$ in $\theta$ and see how $s$ acts on those states.
By the same token, for collapsed transformations, we go through the preimages of
$t$ in $\varphi$, and construct a separate cascade transformations for each.
Thus, all the different actions in $\varphi^{-1}(t)$ contribute to $U$.
For the complete proof, see Theorem \ref{thmcoveringlemma} in Appendix \ref{app:proofs}.

\subsection{From Emulation to Interpretation}
\label{sect:I}
When working with decompositions,  we often want to find the meaning of some coordinatized computations, i.e., doing the computation in the cascade product and see what it corresponds to in the original semigroup.
We call this use of a surjective morphism an \emph{interpretation}.

Once we have an emulation, then we get the interpretation for free, since we can always reverse a morphism.
Reversing an injective one (the emulation) gives us a surjective one, the interpretation:
$$\text{from } \ (X,S)\xhookrightarrow{E(\psi,\mu)}(Y,T)\cp(Z,U)\ \text{ to } \ (X,S)\xtwoheadleftarrow{I(\psi^{-1},\mu^{-1})}(Y,T)\cp(Z,U).$$
If we have $E$ computed completely, then it is possibly to invert the relation.
However, it is better to compute inverse relations directly:
\begin{align*}
\psi^{-1}(y,z)&= z\cdot w_{y}^{-1},\\
\mu^{-1}(t, d_{s,t})&=\{ w_y u w_{yt}^{-1}\mid y\in Y' .\subseteq\image\theta\}
\end{align*}
We need only a subset $Y'$ of $Y$ such that $\bigcup_{y\in Y'}\theta^{-1}(y)=X$,
since the actions are the same on the same points.
The fact that $EI$ should be the identity on $(X,S)$ provides a tool for the
verification of a computational implementation (for details see Corollary
\ref{cor:EI} in Appendix \ref{app:proofs}).

\subsection{Practicality -- Working with the Generators}

A degree $n$ transformation semigroup can have up to $n^n$, thus an efficient algorithm should work on a small generator set.
Here, we check that if $R(\theta,\varphi)$ is defined for $(X,S=\langle
A\rangle)$ only on the generators, we can compute $E(\psi,\mu)$ by lifting those generators only.

The relation on states $\theta$ and its inverse $\theta^{-1}$ should be fully represented, i.e., all $\theta(x)$ values known and stored explicitly.
For computing $\psi(x)$, we need $\theta^{-1}$ for the labelling.

For lifting a single transformation $s$, for computing $\mu(s)$, we need in
addition $\varphi(s)$ computed.
However, we do not need any other values of $\varphi$.
Thus, we can lift each generator without fully computing $\varphi$.
Moreover, the same applies to $I(\psi^{-1},\mu^{-1})$, thus only the relations on states need to be fully enumerated.

\subsection{What is $U$ exactly?}
\label{sect:Uy}
In general, the above construction does not yield a semigroup for $U$, but a \emph{semigroupoid} with partially defined multiplication.
The transformations $u:Z\mapsto Z$ are technically composable, but if we naively compute what they generate we may get a semigroup bigger than what is needed for the emulation.
If the labelling functions do not cancel out, they introduce their arbitrary mappings (unrelated to $S$) into $U$, violating efficiency and correctness.

We can identify the `islands of composability' in $U$ and relate them back to $S$.
The compatible elements are defined by the key idea  that \emph{labelling functions should cancel out in the middle.}
We have a family of transformation semigroups $(Z,U_y)$ for all $y\in \image\theta$, where $U_y=\{w_y^{-1}sw_y\mid yt=y, s\in\varphi^{-1}(t)\}$.
For instance, let $s=s_1s_2\cdots s_k$ be expressed as a product of generators, and $t_1t_2\cdots t_k$ be a corresponding product, such that $t_i\in\varphi(s_i)$.
If $y\cdot t_1t_2\cdots t_k=y$, then in the bottom component we may do context switchings, but they all cancel out: $w_y^{-1}s_1w_{yt_1}\cdot w_{yt_1}^{-1}s_2 w_{yt_1t_2}\cdots w^{-1}_{t_1t_2\cdots t_{k-1}}s_kw_{y}=w_y^{-1}s_1s_2\cdots s_kw_y\in U_y$.

To show that $U_y$ does not contain any maps not present in $S$, we define $f_y: U_y \hookrightarrow S$ by $f_y(u)=w_yuw_y^{-1}$.
 Thus, if $u=w_y^{-1}sw_y$, then $f(u)=w_yw_y^{-1}sw_yw_y^{-1}=s$.
It is easy to check that $f_y(u_1)\cdot f_y(u_2)=f_y(u_1\cdot u_2)$.
Therefore, we can build an injective relational morphism $V_y(w_y^{-1}, f_y):(Z,U_y)\hookrightarrow(X,S)$ to verify that $U_y$ embeds into $S$.
This satisfies the symmetry principle of a good decomposition: the cascade product emulates, while its components are emulated by the original semigroup.

We can characterize the images of $f_y$ as the stabilizers of $\theta^{-1}(y)$ in $S$, since $yt=y$ implies $\theta^{-1}(y)\varphi^{-1}(t)\subseteq \theta^{-1}(yt)=\theta^{-1}(y)$.
This is how local computation is expressed in the bottom level component.

\section{The Specifics}
\label{sec:specifics}
Here we describe how to create surjective morphisms, the input for the
algorithm, and give labelling functions.
These specific details are parameters of the main algorithm.

\subsection{The $n(n-1)$ Method}
\label{sect:ntimesnminusone}
We map each state $X=\{1,2,\ldots,n\}$ to a set of states by
$\theta: X\rightarrow 2^X$.
The default choice in \cite{FromRelToEmulationNehaniv1996} is
$\theta(x)=X\setminus\{x\}$.
In other words, each state goes to the set of all states but itself.
Note that these sets of lifts are overlapping and $\theta$ has the property of being a self-inverse, i.e., $\theta=\theta^{-1}$.
For transformations, we use
\begin{equation*}
    \varphi(s)=
        \begin{cases}
            \{s\} & \text{if } X\cdot s=X,  s\text{ is a permutation,}\\
            \{c_j \mid j\notin \image s\} & \text{if }X\cdot s\subsetneq X, s \text{ collapses some states.}
        \end{cases}
\end{equation*}
These `backwards' relations have a simple explanation. By composing functions we can only reduce the size of the image. To make the morphism work, we need a bigger $\varphi(s)$ set for an $s$ with a smaller image set. 

Now we have $R(\theta, \varphi)$ defining a surjective morphism, but we still
need a corresponding labelling method.
Since $\theta$ misses only a single point, the labelling of states only needs to
deal with just a single hole in the state set.
We need to label $n-1$ states.
For a top level state $y$ we define the map
$w_y:\theta^{-1}(y) \to Z$ by
\begin{equation*}
  w_y(x)=\begin{cases}
    x & \text{if } x < y,\\
    x-1  & \text{if } x > y.
  \end{cases}
\end{equation*}
This gives a $n(n-1)$ cascade transformation representation of a degree $n$
transformation semigroup, since $|\theta^{-1}(y)|=n-1$ for all $y$.
The canonical example for this method is the decomposition of the full transformation semigroup.
For arbitrary transformation semigroups, it does not produce
an efficient separation of computation.
Transformations may get repeated on the lower level.

A simple example can demonstrate this problem. Consider
the semigroup generated by the transformation
$\big(\begin{smallmatrix}
  1 & 2 & 3 & 4 & 5\\
  2 & 3 & 1 & 5 & 4
\end{smallmatrix}\big)$, which is the permutation $(1,2,3)(4,5)$. Thus, we have a
group, and all the elements will make it to the top level.
At the same time, the 2- and the 3-cycles will also make it to the bottom level,
since the labelling has a cut-off at four states: $|\theta^{-1}|=4$.
Therefore, we need to find a better method.

\subsection{Generalized Labelling: The Squashing Function}
\label{sect:squashing}
For a general surjective morphism $R(\theta, \varphi)$, the preimages of
$\theta$ can be of different size, thus simply removing a hole is not enough.
There may be more missing points in the preimage set, so need to `squash' the set
to remove all the holes.

Let $\theta^{-1}(y)=\{x_1, x_2, \ldots, x_k\}$.
Let $\sigma$ be the permutation that orders these states, so $x_{\sigma(1)} < x_{\sigma(2)} < \ldots < x_{\sigma(k)}$.
Then $w_y=\sigma^{-1}$.
The number of states used in $Z$ is the size of the preimage.
If we need $k$ states, then we use the first $k$ states in Z.
Consequently, $|Z|=\max_y|\theta^{-1}(y)|$.
For instance, if $\theta^{-1}(y)=\{4,5,2\}$ and $ \theta^{-1}(yt)=\{5,1\}$, then
 three states would be sufficient in $Z$ for the squashing function to operate.
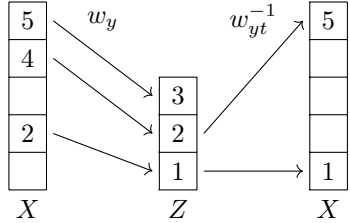
\begin{center}
\begin{tikzpicture}[inner sep=7pt]
\draw[step=0.5] (0,0) grid (0.5,2.5);
\node (a2) at (0.25,0.75) {2};
\node (a4) at (0.25,1.75) {4};
\node (a5) at (0.25,2.25) {5};
\node at (1.25,2.25) {$w_y$};
\node at (0.25,-0.25) {$X$};

\draw[step=0.5] (1.99,0) grid (2.5,1.5); 
\node (b1) at (2.25,0.25) {1};
\node (b2) at (2.25,0.75) {2};
\node (b3) at (2.25,1.25) {3};
\node at (2.25,-0.25) {$Z$};
\node at (3.25,2.25) {$w_{yt}^{-1}$};

\draw[step=0.5] (3.99,0) grid (4.5,2.5);
\node (c1) at (4.25,0.25) {1};
\node (c5) at (4.25,2.25) {5};
\node at (4.25,-0.25) {$X$};

\draw [->] (a2.east) -- (b1.west);
\draw [->] (a4.east) -- (b2.west);
\draw [->] (a5.east) -- (b3.west);
\draw [->] (b1.east) -- (c1.west);
\draw [->] (b2.east) -- (c5.west);
\end{tikzpicture}
\end{center}

\subsection{Getting a surjective morphism}
Arguably, the Covering Lemma only does half of the decomposition.
It needs an
input surjective morphism for defining the top level component.
Constructing this morphism could already contribute to the understanding of
$(X,S)$, and in certain applications it might be defined naturally.
The above $n(n-1)$ method gives a working, but not necessarily balanced
  decomposition (see \ref{sect:ntimesnminusone}) and we can also project a
  holonomy decomposition into its top level \cite{holonomy_algorithm}.
  Here we describe two more methods: congruences and local monoids.

An equivalence relation $\sim$ on $X$ is a \emph{right congruence} if $x_1\sim
x_2\implies x_1s\sim x_2s$ for all $s\in S$.
In other words, the semigroup action is compatible with the equivalence classes.
Thus, we have $(X/{\sim}, S')$ well-defined, where $S'$ is $S$ made faithful,
and the morphism to this is injective.
In practice, we can find a congruence by a standard method.
From the application domain we may have some intuition about which states to
identify.
We can specify these in disjoint sets of states, and all the remaining states
form singleton classes.
Then, an algorithm dual to the classical DFA minimization algorithm can compute
the minimal (finest) congruence containing the input identified states in
equivalence classes, acting as a closure  operator.
We check the action of the generators for compatibility with the existing
classes.
If elements from a class go to different classes, then we need to merge those.
We iterate this until there is no merging.
In general, there are many such congruences for a transformation semigroup, but in special cases we may have only
the trivial congruence ($X$ itself is the only class), as for the full transformation semigroup.

The \emph{local transformation monoid} $(Xe, eSe)$ is defined for an idempotent $e=e^2\in
S$.
It localizes the semigroup's action to a subset $Xe$ of the state set.
It is an excellent tool for finding an interesting subsemigroup by mapping $s\mapsto ese$, but
that is not necessarily a homomorphic image: $st\mapsto este$, and $ese$ and $ete$ have product $eseets=esete$, but
$set\neq st$ in general. This method works for commutative semigroups and in
some special cases.

\section{Examples}

\subsection{Simple cycle collapsing}
This example provides a minimal but still meaningful demonstration of the action transfer algorithm.
We only have three states, thus we can denote the transformation by their image list in condensed format.
For instance, $\big(\begin{smallmatrix}
  1 & 2 & 3\\
  1 & 3 & 2
\end{smallmatrix}\big)$
will be written as 132.
With this notation we define $(X,S)$ as $X=\{1,2,3\}$ and $S=\{e=123, p=132,
  c_1=111, c_2=222, c_3=333\}$. The element $p$ has a non-trivial permutation.
The target $(Y,T)$ is defined by $Y=\{1,2\}$ and $T=\{e'=12, c_1'=11, c_2'=22\}$.
Now $\theta$ collapses 2 and 3, and fixes 1, i.e., $\theta(2)=\theta(3)=\{2\}$, $\theta(1)=\{1\}$.
This determines $\varphi$ as well, e.g., $\varphi(c_2)=\varphi(c_3)=\{c_2'\}$, $\varphi(p)=\{e'\}$.
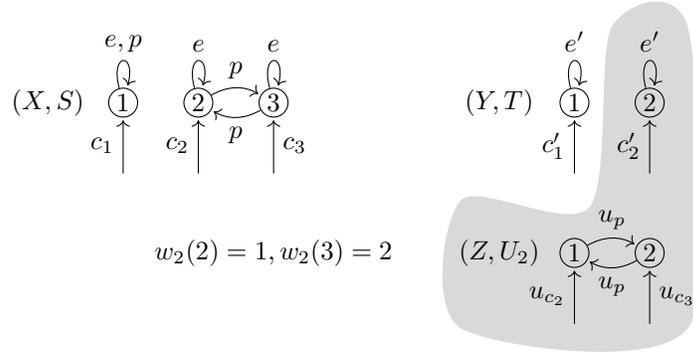
\begin{center}
\usetikzlibrary {automata,positioning}
\begin{tikzpicture}[shorten >=1pt,node distance=1cm,
    every state/.style={inner sep=1pt,minimum size=2pt},on grid,auto,
    ]
    \filldraw [gray!30] plot [smooth cycle] coordinates {(6.5,1)(7.5,1)(7.5,-3)(4.5,-3)(4.5,-1.3)(6.2,-1.4)};
  \node[state]  (q_1)                      {1};
  \node (XS) [left=of q_1] {$(X,S)$};
  \node[state]          (q_2) [right=of q_1] {2};
  \node[state]          (q_3) [right=of q_2] {3};
  \node[state] (t_1) at (6,0) {1};
  \node (YT) [left=of t_1] {$(Y,T)$};
  \node[state] (t_2) [right=of t_1] {2};
  \node[state,draw=none] (c_1) [below=of q_1] {};
  \node[state,draw=none] (c_2) [below=of q_2] {};
  \node[state,draw=none] (c_3) [below=of q_3] {};
  \node[state,draw=none] (tc_1) [below=of t_1] {};
  \node[state,draw=none] (tc_2) [below=of t_2] {};
  \node[state] (z_1) [below=of tc_1] {1};
   \node[state] (z_2) [below=of tc_2] {2};
   \node[state,draw=none] (zc_1) [below=of z_1] {};
   \node[state,draw=none] (zc_2) [below=of z_2] {};
   \node (ZU) [left=of z_1] {$(Z,U_2)$};
   \node (wu) [below=of c_3] {$w_2(2)=1, w_2(3)=2$};
\path[->] (q_1) edge [loop above]  node        {$e,p$} ()
  (q_2) edge [bend left] node {$p$}(q_3)
  (q_3) edge [bend left] node {$p$} (q_2)
  (q_2) edge [loop above] node {$e$} ()
  (q_3) edge [loop above] node {$e$} ()
  (t_1) edge [loop above] node {$e'$}()
  (t_2) edge [loop above] node {$e'$} ()
  (tc_1) edge node {$c'_1$} (t_1)
  (tc_2) edge node {$c'_2$} (t_2)
  (zc_1) edge node {$u_{c_2}$} (z_1)
  (zc_2) edge node [right] {$u_{c_3}$} (z_2)
  (z_2) edge [bend left] node {$u_p$}(z_1)
  (z_1) edge [bend left] node {$u_p$} (z_2)

  (c_1) edge node {$c_1$} (q_1)
  (c_2) edge node {$c_2$} (q_2)
  (c_3) edge node [right] {$c_3$}  (q_3);
\end{tikzpicture}
\end{center}
Note that $R(\theta,\varphi)$ is a surjective homomorphism since all the images
are singletons.
What is the information lost? All the original transformations that had
different actions on $\{2,3\}$ cannot be distinguished any more.
Most notably, we should get a permutation on the bottom level.
Indeed, when lifting $p$, we have a non-trivial dependency from $y=2$ to the
permutation $(1,2)$, as the original cycle $(2,3)$ is relabelled.
Since $|\theta^{-1}(1)|=1$, there are no fine-grained details attached to this state, $U_1$ is just the trivial monoid.

\subsection{``Bad'' examples}
\label{section:bad}
We consider a resulting decomposition ``bad'' if it does not do a balanced
distribution of computation over the components.
For the usability of the decomposition we need to follow the `Divide and conquer!' principle for dealing with complexity.
In the extreme cases the decomposition only adds the overhead of the hierarchical form.

For a surjective homomorphism that collapses everything down to the trivial
monoid, i.e., $\forall x\in X, \theta(x)=\{1\}$ and $\forall s\in S, \varphi(s)=\{e\}$ ($e$ is the identity in $T$), we get everything on the second level through a dummy dependency
function with a single input value, since $\theta^{-1}(1)=X$.
As the other extreme, if we have an isomorphism to start with, we will get a trivial bottom level component (see Fact \ref{fact:blockedaction}).

\section{Comparison: Holonomy versus Covering Lemma}
The decisive difference is that holonomy is a monolithic method, while the Covering Lemma allows incremental decompositions.
Both work with generators, but the holonomy decomposition constructs a detailed analysis of the set of images of the semigroup action, requiring to compute possibly as many as $2^n$ subsets of the state set.

The holonomy algorithm does provide more information.
It gives compression: identifying equivalence classes of subsets of the state set with isomorphic
groups action on them the mappings within the classes.
When using the covering lemma method, this information needs to be computed separately for the local transformation semigroups $\{U_y\mid y\in Y\}$.
However, this is easier than in holonomy, since we can use the state set $Y$ for finding the strongly connected components in reachability relation under $T$, instead of working in the set of image sets.

Another seemingly big difference is that in the covering lemma method we do not separate
the parallel components, as they are defined on the same set of states.
The holonomy uses disjoint union, so we can have many state components on one level.
However, this can be easily replicated, if needed, in the covering lemma method.
We could modify the labelling function not to reuse the states.

\section{Conclusions}
Here we extended the computer algebra toolkit for the hierarchical decompositions
of finite transformation semigroups by a simplified kernel-based two-level algorithm.
In order to meet the requirements of the applications (as opposed to, or in addition to proving the lemma), we needed to make a few innovations:
generalized labelling for states (the squashing function, \ref{sect:squashing});
explicit computation of the reversed emulation, the \emph{interpretation} (\ref{sect:I}), for typical use cases in applications and for verification purposes;
identifying the family of the bottom level transformation semigroups (\ref{sect:Uy}) in lieu of the derived semigroupoid; 
separation of the general mechanism (Section \ref{sect:generalalg}) from the
specific choices (Section \ref{sec:specifics}).

The holonomy decomposition first implemented in \texttt{SgpDec}
\cite{SgpDec2014} allows one to make the transition from having only hand-calculated small decomposition examples to the machine computed and verified coordinatizations of real-world dynamical systems \cite{PhilTransA2015}.
Similarly, the now implemented Covering Lemma algorithm (see Appendix \ref{app:software}) makes the transition from the fixed highest resolution decompositions only, monolithic calculation tool to a more flexible one.
By iterating decompositions it extends the `computational horizon' as we can compute and study bigger examples.
By allowing an arbitrary surjective morphism as its input, it enables `precision engineering', constructing tailor-made decompositions for particular discrete dynamical systems.
We made effort to present the algorithm in an elementary way accessible to a wider audience.
We hope this will allow other researchers and students to apply the
decompositions in their projects.

%
%

\bibliographystyle{plain}

\bibliography{coords}
\appendix

\section{Proofs}
\label{app:proofs}
Here we state the decomposition part of Covering Lemma without the derived category.
The main idea
is the context-dependent partial transfer of the original semigroup action, thus
the proof merely checks that the resulting constructions satisfies the definition of emulation.
\begin{theorem}[Covering Lemma]\label{thmcoveringlemma}
If $R(\theta,\varphi):(X,S)\twoheadrightarrow (Y,T)$ is a surjective relational morphism,
then there exists an emulation, an injective relational morphism $E(\psi,\mu): (X,S)\hookrightarrow (Y,T)\cp(Z,U)$.
\end{theorem}

\proof
Let's choose a state $x\in X$ and a transformation $s\in S$.
The lifted states are coordinate pairs $\psi(x)=\{ (y,xw_y) \mid y \in
\theta(x)\}$.
The lifted transformations are transformation cascades
$\mu(s)=\{(t, w_y^{-1} s w_{yt}) \mid t\in \varphi(s)\}$,
where $w_y:X\to Z$ is any partial one-to-one function defined on $\theta^{-1}(y)$ for all $y\in\image\theta$.
Then $\psi(x)\cdot\mu(s)=\{(yt, xw_y\cdot  w_y^{-1} s w_{yt})\mid
y\in\theta(x),t\in\varphi(s)\}$ according cascade multiplication.
Thus, we have the set of coordinate pairs $(yt, xsw_{yt})$ ranging through all the lifts of
$x$ and $s$.

On the other hand, $\psi(x\cdot s)=\{(y',xsw_{y'})\mid y'\in\theta(xs)\}$.
Since $R(\theta,\varphi)$ is a relational morphism, we have $\theta(x)\cdot
\varphi(s)\subseteq \theta(x\cdot s)$.
Consequently, the set of $yt$ top level states above is a subset of the set of these $y'$ states.
Therefore, $\psi(x)\cdot\mu(s)\subseteq \psi(x\cdot s)$, and $E(\psi,\mu)$ is a relational morphism.

For injectivity, pick two states $x_1,x_2\in X$.
Assume that injectivity does not hold for $\psi$ and we have a coordinate pair
such that $(y,z)\in
\psi(x_1)\cap\psi(x_2)$.
Now, $z$ can arise in two ways, $x_1\mapsto (y, x_1w_y)$ and $x_2 \mapsto (y, x_2w_y)$.
Since $w_y$ is one-to-one, $x_1w_y=w_2w_y\implies x_1=x_2$.
Therefore, $E$ is injective on states.

Let's pick two transformations, $s_1, s_2$.
Assume that $\mu(s_1)\cap \mu(s_2)\neq\varnothing$.
Then choose a cascade transformation $(t,w_y^{-1}sw_{yt})$ from the intersection.
This again can arise in two ways: $s_1\mapsto w_y^{-1}s_1w_{yt}$ and $s_2\mapsto
w_y^{-1}s_2w_{yt}$.
Since these are just different expressions of the same cascade transformation,
$\forall y\in\image\theta$ we have $w_y^{-1}s_1w_{yt} = w_y^{-1}s_2w_{yt}$.
This implies $s_1=s_2$ by faithfulness of $(X,S)$, therefore $E$ is injective on transformations.
\qed

\begin{corollary}[Identities from Composite Relations]\label{cor:EI}
Inverting the emulation $E(\psi, \mu)$ gives the surjective relational morphism
$I(\psi^{-1},\mu^{-1}): (Y,T)\cp(Z,U)\twoheadrightarrow (X,S)$, and $IE$ is the
identity on $(X,S)$, and $EI$ is the identity on the cascade product.
\end{corollary}

\proof
It is a standard property of relational morphisms that the inverse of an injective
morphism is a surjective one.
It is also elementary that a relation combined with its inverse gives the
identity relation.
However, it is instructional to check this with the specific formulas.
For states,
$IE(x)=I(E(x))=\psi^{-1}(\psi(x))=\psi^{-1}((y,xw_y))=xw_yw_y^{-1}=x, \forall
y\in\theta(x)$, and
$EI((y,z))=E(I((y,z)))=\psi(\psi^{-1}((y,z)))=\psi(zw_y^{-1})=(y, zw_y^{-1}w_y)=(y,z)$.
For transformations,
$IE(s)=I(E(s))=\mu^{-1}(\mu(s))=\mu^{-1}((t,w_y^{-1}sw_{yt}))=w_yw_y^{-1}sw_{yt}w_{yt}^{-1}=s$,
and  using the local action $u\in U$ we have
$EI((t,u))=E(I((t,u)))=\mu(\mu^{-1}((t,u)))=\mu(w_yuw_{yt}^{-1})=(t, w_y^{-1}w_yuw_{yt}w_{yt}^{-1})=(t,u)$.
\qed

\vspace{1em}
The set of preimages $\theta^{-1}(y)$ serve as the channel for transferring the
action in $(X,S)$ to $(Z,U)$.
If $|\theta^{-1}(y)|=1$ for all $y\in Y$, then we completely block the way for any
nontrivial action to get in the bottom component.
This may seem like a way to craft a counterexample to the Covering Lemma, but
here we show that in that case the top level component can emulate $(X,S)$ in
itself.
\begin{fact}[Blocked action transfer]\label{fact:blockedaction}
If $\theta$ is bijective (as a function), i.e., $\theta(x)=\{y\}$ is a singleton
$\forall x\in X$, and  $\theta(x_1)=\theta(x_2)\implies x_1=x_2$, so $|X|=|Y|$,
then $\varphi$ is injective.
\end{fact}
\proof
Assume that $t\in\varphi(s_1)\cap\varphi(s_2)$. Since $\theta^{-1}$ is also
bijective, the action of $t$ on $Y$ uniquely determines the action of $X$
: for all $y\in Y$, for $yt$ we have
$\theta^{-1}(y)s_1=\theta^{-1}(y)s_2$, thus $s_1=s_2$ by faithfulness of $(X,S)$.\qed

\section{Software implementation}
\label{app:software}
The implementation of the algorithms presented here are available\footnote{At
  time of writing this early draft it may still be in the experimental
  public repository \url{https://github.com/egri-nagy/sgplab} under the
  \texttt{CoveringLemma} folder.} in the
\texttt{SgpDec} package \cite{SgpDec} written for the \texttt{GAP} \cite{GAP4}
computer algebra system relying on the foundational \texttt{Semigroups} package
\cite{Semigroups}.
To implement relational morphisms we use the hash-maps provided by the
\texttt{datastructures} \cite{datastructures} package.
This package has a
lightweight reference implementation for hash-maps and comes with a recursive hash code computation for
composite data structures, thus transformation cascades in \texttt{SgpDec} can
be used as keys.
We also use hash-maps when implementing the partial labelling functions and the congruence closure algorithm.

For representing transformations in the local $U_y$ components, we use identities for any undefined maps, instead of implementing them as partial transformations  with a sink state.
Both methods can introduce unrelated transformations if we multiply incompatible elements.
The sink state representation has a second disadvantage: it increases the number of states needed.

The source code closely follows this paper.
Function and variable names are as close as possible to the notation in the
mathematical text.
There are several implementations for $\theta$ with their function names starting
with the prefix \texttt{ThetaFor}.
For the $n(n-1)$ method there is \texttt{ThetaForPermutationResets}, for the
congruences \texttt{ThetaForCongruence}, for local monoids
\texttt{ThetaForLocalMonoid}. In case users start with a holonomy decomposition, they
can use \texttt{ThetaForHolonomy}. For testing purposes we provide
\texttt{ThetaForConstant}, mapping everything to a single state.
For all these, we have the matching algorithms for $\varphi$, prefixed by
\texttt{PhiFor}.
Once \texttt{theta} and \texttt{phi} are defined, for computing $\psi$ and $\mu$
the user simply calls \texttt{Psi(theta)} and \texttt{Mu(theta,phi)}.

Here is a sample session to demonstrate the workflow of investigating of an
unknown transformation semigroup (a randomly generated one in this example).
\begin{verbatim}
S:=Semigroup([Transformation([1,6,11,12,11,10,7,13,7,1,2,1,1]),
              Transformation([2,10,3,3,8,7,2,4,5,6,5,3,4])]);
\end{verbatim}
The semigroup has 9221 elements.
We apply the congruence based method for finding a surjective homomorphism onto
a smaller transformation semigroup.
Arbitrarily, we try to identify 1 with 2, and 3 with 4.
\begin{verbatim}
gap> partition := StateSetCongruence(Generators(S), [[1,2],[3,4]]);
[[1,2,6,7,10],[3,4,5,8],[9],[ 11,12,13]]
\end{verbatim}
Now, we know that the homomorphic image will be acting on 4 points. How big is
this semigroup?
\begin{verbatim}
theta := ThetaForCongruence(partition);
phi := PhiForCongruence(partition, Generators(S));
\end{verbatim}
We compute $\varphi$ only for the generators.
\begin{verbatim}
gap> ImageOfHashMapRelation(phi);
[ Transformation( [ 1, 2, 2, 2 ] ), Transformation( [ 1, 4, 1, 1 ] ) ]
\end{verbatim}
This has 5 elements and it is an aperiodic semigroup, i.e., it does not have a non-trivial subgroup.
Any non-trivial group action should appear on the second level of the decomposition.
The generators of the cascade product are given by $\mu$.
\begin{verbatim}
mu := Mu(theta, phi);  ImageOfHashMapRelation(mu);
[ <trans cascade with 2 levels with (4, 5) pts, 5 dependencies>,
  <trans cascade with 2 levels with (4, 5) pts, 5 dependencies> ]
\end{verbatim}
These generators are the following cascade transformations.
\begin{verbatim}
Dependency function of depth 1 with 1 dependencies.
[  ] -> Transformation( [ 1, 2, 2, 2 ] )
Dependency function of depth 2 with 4 dependencies.
[ 1 ] -> Transformation( [ 2, 5, 4, 2, 3 ] )
[ 2 ] -> Transformation( [ 1, 1, 4, 2 ] )
[ 3 ] -> Transformation( [ 3, 2, 3 ] )
[ 4 ] -> Transformation( [ 3, 1, 2 ] )
\end{verbatim}

\begin{verbatim}
Dependency function of depth 1 with 1 dependencies.
[  ] -> Transformation( [ 1, 4, 1, 1 ] )
Dependency function of depth 2 with 4 dependencies.
[ 1 ] -> Transformation( [ 1, 3, 5, 4, 1 ] )
[ 2 ] -> Transformation( [ 1, 2, 1, 3 ] )
[ 3 ] -> Transformation( [ 4, 2, 3, 4 ] )
[ 4 ] -> Transformation( [ 2, 1, 1 ] )
\end{verbatim}
The top level state 3 corresponds to the single original state 9, thus
$|\theta^{-1}(3)|=1$.
Still, a non-identity transformation can appear on the second level (only a
single map, $1\mapsto 3$ in the first, and $1\mapsto 4$ in the second generator), since other
top level states can have bigger preimage sets in $\theta$.

\end{document}